\documentclass[12pt]{article}
\usepackage{fullpage}
\usepackage{amsmath,amsfonts,epsfig,amssymb,amsthm,color,graphicx}
 \theoremstyle{plain}
   \newtheorem{theorem}{Theorem}[section]
   \newtheorem{proposition}[theorem]{Proposition}
   
   \newtheorem{corollary}[theorem]{Corollary}
\theoremstyle{definition}
   \newtheorem{definition}[theorem]{Definition}

\title{The alternating sign matrix polytope}
\author{Jessica Striker\\
\small School of Mathematics\\[-0.8ex]
\small University of Minnesota\\[-0.8ex]
\small Minneapolis, MN 55455\\[-0.8ex]
\small \texttt{jessica@math.umn.edu}}

\begin{document}

\maketitle

\begin{abstract}
We define the alternating sign matrix polytope as the convex hull of $n\times n$ alternating sign matrices and prove its equivalent description in terms of inequalities. This is analogous to the well known result of Birkhoff and von Neumann that the convex hull of the permutation matrices equals the set of all nonnegative doubly stochastic matrices. We count the facets and vertices of the alternating sign matrix polytope and describe its projection to the permutohedron as well as give a complete characterization of its face lattice in terms of modified square ice configurations. Furthermore we prove that the dimension of any face can be easily determined from this characterization.
\end{abstract}

\section{Introduction and background}
\label{sec:back}
The Birkhoff polytope, which we will denote as $B_n$, has been extensively studied and generalized. It is defined as the convex hull of the $n \times n$ permutation matrices as vectors in $\mathbb{R}^{n^2}$. Many analogous polytopes have been studied which are subsets of $B_n$ (see e.g.~\cite{BRUALDI}). In contrast, we study a polytope containing $B_n$. We begin with the following definitions.

\begin{definition}
\label{def:asm}
Alternating sign matrices (ASMs) are square matrices with the following properties:
\begin{itemize}
\item entries $\in \{0,1,-1\}$
\item the entries in each row and column sum to 1
\item nonzero entries in each row and column alternate in sign
\end{itemize}
\end{definition}

\begin{figure}[htbp]
\[
\left( 
\begin{array}{rrr}
1 & 0 & 0 \\
0 & 1 & 0\\
0 & 0 & 1
\end{array} \right)
\left( 
\begin{array}{rrr}
1 & 0 & 0 \\
0 & 0 & 1\\
0 & 1 & 0
\end{array} \right)
\left( 
\begin{array}{rrr}
0 & 1 & 0 \\
1 & 0 & 0\\
0 & 0 & 1
\end{array} \right)
\left( 
\begin{array}{rrr}
0 & 1 & 0 \\
1 & -1 & 1\\
0 & 1 & 0
\end{array} \right)
\]

\[
\left( 
\begin{array}{rrr}
0 & 1 & 0 \\
0 & 0 & 1\\
1 & 0 & 0
\end{array} \right)
\left( 
\begin{array}{rrr}
0 & 0 & 1 \\
1 & 0 & 0\\
0 & 1 & 0
\end{array} \right)
\left( 
\begin{array}{rrr}
0 & 0 & 1 \\
0 & 1 & 0\\
1 & 0 & 0
\end{array} \right)
\]
\label{fig:n3asm}
\caption{The $3\times 3$ ASMs}
\end{figure}

The total number of $n \times n$ alternating sign matrices is given by the expression
\begin{equation}
\label{eq:product}
\displaystyle\prod_{j=0}^{n-1} \frac{(3j+1)!}{(n+j)!}.
\end{equation}
Mills, Robbins, and Rumsey conjectured this formula~\cite{MRRASMDPP}, and then over a decade later Doron Zeilberger proved it~\cite{ZEILASM}.
Shortly thereafter, Kuperberg found a bijection between ASMs and the statistical physics model of square ice with domain wall boundary conditions (which is very similar to the \emph{simple flow grids} defined later in this paper), and gave a shorter proof using insights from physics~\cite{KUP_ASM_CONJ}. For a detailed exposition of the conjecture and proof of the enumeration of ASMs, see~\cite{BRESSOUDBOOK}.
See Figure~\ref{fig:n3asm} for the seven $3\times 3$ ASMs.  

\begin{definition}
The $n$th alternating sign matrix polytope, which we will denote as $ASM_n$, is the convex hull in $\mathbb{R}^{n^2}$ of the $n\times n$ alternating sign matrices.
\end{definition}

From Definition~\ref{def:asm} we see that permutation matrices are the alternating sign matrices whose entries are nonnegative. Thus $B_n$ is contained in $ASM_n$. The connection between permutation matrices and ASMs is much deeper than simply containment. There exists a partial ordering on alternating sign matrices that is a distributive lattice. This lattice contains as a subposet the Bruhat order on the symmetric group, and in fact, it is the smallest lattice that does so (i.e.\ it is the MacNeille completion of the Bruhat order)~\cite{TREILLIS}. Given this close relationship between permutations and ASMs it is natural to hope for theorems for $ASM_n$ analogous to those known for $B_n$. 

In this paper we find analogues for $ASM_n$ of the following theorems about the Birkhoff polytope (see the discussion in~\cite{BRUALDI} and~\cite{YEMELICHEV}). 

\begin{itemize}
\item $B_n$ consists of the $n\times n$ nonnegative doubly stochastic matrices (square matrices with nonnegative real entries whose rows and columns sum to 1). 
\item The dimension of $B_n$ is $(n-1)^2$. 
\item $B_n$ has $n!$ vertices.
\item $B_n$ has $n^2$ facets (for $n\ge3$) where each facet is made up of all nonnegative doubly stochastic matrices with a 0 
in a specified entry. 
\item $B_n$ projects onto the permutohedron.
\item There exists a nice characterization of its face lattice in terms of elementary bipartite graphs~\cite{BILLERA}.
\end{itemize}

As we shall see in Theorem~\ref{thm:asmchasing}, the row and column sums of every matrix in $ASM_n$ must equal 1. Thus the dimension of $ASM_n$ is $(n-1)^2$ because, just as for the Birkhoff polytope, the last entry in each row and column is determined to be precisely what is needed to make that row or column sum equal 1. 
In Section~\ref{sec:properties} we prove that $ASM_n$ has $4[(n-2)^2+1]$ facets and its 
vertices are the alternating sign matrices.
We also prove analogous theorems about the inequality description of $ASM_n$ (Section~\ref{sec:asm}), the face lattice (Section~\ref{sec:facelattice}), and the projection to the permutohedron (Section~\ref{sec:properties}). See~\cite{ZIEGLER} for background and terminology on polytopes.

The alternating sign matrix polytope was independently defined in~\cite{KNIGHT} in which the authors also study the integer points in the $r$th dilate of $ASM_n$ calling them \emph{higher spin alternating sign matrices}.

\section{The inequality description of the ASM polytope}
\label{sec:asm}

The main theorem about the Birkhoff polytope is the theorem of Birkhoff~\cite{BIRKHOFFPOLY} and von~Neumann~\cite{VONNEUMANN}
 which says that the Birkhoff polytope can be described not only as the convex hull of the permutation matrices but equivalently  
as the set of all nonnegative doubly stochastic matrices (real square matrices with row and column sums equaling 1 whose entries are nonnegative). 

The inequality description of the alternating sign matrix polytope is similar to that of the Birkhoff polytope. 
It consists of the subset of doubly stochastic matrices 
(now allowed to have negative entries) whose partial sums in each row and 
column are between 0 and 1. The proof uses the idea of von Neumann's proof
of the inequality description of the Birkhoff polytope~\cite{VONNEUMANN}. 

Note that in~\cite{KNIGHT} Behrend and Knight approach the equivalence of the convex hull definition and the inequality description of $ASM_n$ in the opposite manner, defining the alternating sign matrix polytope in terms of inequalities and then proving that the vertices are the alternating sign matrices. 

\begin{theorem}
\label{thm:asmchasing}

The convex hull of $n\times n$ alternating sign matrices consists of all 
$n\times n$ real matrices $X=\{x_{ij}\}$ such that:

\begin{align}
\label{eq:partialcolumnsum}
0 \le \sum_{i=1}^{i'} x_{ij}&\le 1 \hspace{.6in} \forall\mbox{ } 1\le i'\le n, 1\le j\le n.\\
\label{eq:partialrowsum}
0 \le \sum_{j=1}^{j'} x_{ij}&\le 1 \hspace{.6in} \forall\mbox{ } 1\le j'\le n, 1\le i\le n.\displaybreak\\
\label{eq:sumagain}
\sum_{i=1}^n x_{ij}&= 1 \hspace{.6in} \forall\mbox{ }1\le j \le n.\\
\label{eq:sumagain2}
\sum_{j=1}^n x_{ij}&= 1 \hspace{.6in} \forall\mbox{ } 1\le i \le n.
\end{align}
\end{theorem}
\begin{proof}
Call the subset of $\mathbb{R}^{n^2}$ given by the above inequalities $P(n)$.
It is easy to check that the convex hull of the alternating sign 
matrices is contained in the set $P(n)$. 
It remains to show that any $X\in P(n)$ can be
written as a convex combination of alternating sign matrices.

Let $X\in P(n)$. Let $r_{ij}=\sum_{j'=1}^j x_{ij}$ and $c_{ij}=\sum_{i'=1}^i x_{ij}$. Thus the $r_{ij}$ are the row partial sums and the $c_{ij}$ are the column partial sums. It follows from (\ref{eq:partialcolumnsum}) and (\ref{eq:partialrowsum}) that $0\le 
r_{ij},c_{ij}\le1$ for all $1\le i,j\le n$. Also, from (\ref{eq:sumagain}) and (\ref{eq:sumagain2}) we see that $r_{in}=c_{nj}=1$ for all $1\le i,j\le n$. If we set $r_{i0}=c_{0j}=0$ we see that every entry $x_{ij}\in X$ satisfies $x_{ij}=r_{ij}-r_{i,j-1}=c_{ij}-c_{i-1,j}$. Thus 
\begin{equation}
\label{eq:rcrel}
r_{ij}+c_{i-1,j}=c_{ij}+r_{i,j-1}. 
\end{equation}

Using von Neumann's terminology, we call a real number $\alpha$ \emph{inner} if $0<\alpha<1$. We construct a circuit in $X$ such that the partial sum between adjacent matrix entries in the circuit be an inner. So we rewrite the matrix $X$ with the partial sums between entries as shown below.

\[
\left(
\begin{array}{ccccccccccc}
& \textcolor{blue}{c_{01}}&&\textcolor{blue}{c_{02}}&& &&\textcolor{blue}{c_{0,n-1}}&&\textcolor{blue}{c_{0n}}&\\
\textcolor{blue}{r_{10}}&x_{11}&\textcolor{blue}{r_{11}}&x_{12}&\textcolor{blue}{r_{12}}& &&x_{1,n-1}&\textcolor{blue}{r_{1,n-1}}&x_{1n}&\textcolor{blue}{r_{1n}}\\
& \textcolor{blue}{c_{11}}&&\textcolor{blue}{c_{12}}&&\ldots&&\textcolor{blue}{c_{1,n-1}}&&\textcolor{blue}{c_{1n}}&\\
\textcolor{blue}{r_{20}}&x_{21}&\textcolor{blue}{r_{21}}&x_{22}&\textcolor{blue}{r_{22}}& &&x_{2,n-1}&\textcolor{blue}{r_{2,n-1}}&x_{2n}&\textcolor{blue}{r_{2n}}\\
&&&&&&&&&&\\
&&\vdots&&&&&&\vdots&&\\
& \textcolor{blue}{c_{n-1,1}}&&\textcolor{blue}{c_{n-1,2}}&& &&\textcolor{blue}{c_{n-1,n-1}}&&\textcolor{blue}{c_{n-1,n}}&\\
\textcolor{blue}{r_{n0}}&x_{n1}&\textcolor{blue}{r_{n1}}&x_{n2}&\textcolor{blue}{r_{n2}}&\ldots&&x_{n,n-1}&\textcolor{blue}{r_{n,n-1}}&x_{nn}&\textcolor{blue}{r_{nn}}\\
& \textcolor{blue}{c_{n1}}&&\textcolor{blue}{c_{n2}}&&&&\textcolor{blue}{c_{n,n-1}}&&\textcolor{blue}{c_{nn}}&
\end{array}
\right)
\]

Begin at the vertex to the left or above any inner partial sum; if no such partial sum exists, then $X$ is an alternating sign matrix. Then there exists an adjacent inner partial sum by (\ref{eq:rcrel}). By repeated application of (\ref{eq:rcrel}) to each new inner partial sum, a path can then be formed by moving from entry to entry of $X$ along inner partial sums. Since $X$ is of finite size and all the boundary partial sums are 0 or 1 (i.e. non--inner), the path eventually reaches an entry in the same row or column as a previous entry yielding a circuit in $X$ whose partial sums are all inner. Using this circuit we can write $X$ as a convex combination of two matrices in $P(n)$, each with at least one more non--inner partial sum, in the following way.

\begin{figure}
\centering
$\begin{array}{lcr}
\left(
\begin{array}{rrrrr}
0&.4&.5&.1&0\\
.4&-.4&.5&0&.5\\
.6&.4&-.3&-.1&.4\\
0&.3&-.3&.9&.1\\
0&.3&.6&.1&0
\end{array}
\right)
&
\Rightarrow
&
\left(
\begin{array}{ccccccccccc}
& \textcolor{blue}{0}&&\textcolor{blue}{0}&&\textcolor{blue}{0}&&\textcolor{blue}{0}&&\textcolor{blue}{0}&\\
\textcolor{blue}{0}&0&\textcolor{blue}{0}&.4&\textcolor{red}{\textbf{.4}}&.5&\textcolor{red}{\textbf{.9}}&.1&\textcolor{blue}{1}&0&\textcolor{blue}{1}\\
& \textcolor{blue}{0}&&\textcolor{red}{\textbf{.4}}&&\textcolor{blue}{.5}&&\textcolor{red}{\textbf{.1}}&&\textcolor{blue}{0}&\\
\textcolor{blue}{0}&.4&\textcolor{red}{\textbf{.4}}&-.4&\textcolor{blue}{0}&.5&\textcolor{blue}{.5}&0&\textcolor{blue}{.5}&.5&\textcolor{blue}{1}\\
& \textcolor{red}{\textbf{.4}}&&\textcolor{blue}{0}&&\textcolor{blue}{1}&&\textcolor{red}{\textbf{.1}}&&\textcolor{blue}{.5}&\\
\textcolor{blue}{0}&.6&\textcolor{red}{\textbf{.6}}&.4&\textcolor{blue}{1}&-.3&\textcolor{red}{\textbf{.7}}&-.1&\textcolor{blue}{.6}&.4&\textcolor{blue}{1}\\
& \textcolor{blue}{1}&&\textcolor{red}{\textbf{.4}}&&\textcolor{red}{\textbf{.7}}&&\textcolor{blue}{0}&&\textcolor{blue}{.9}&\\
\textcolor{blue}{0}&0&\textcolor{blue}{0}&.3&\textcolor{red}{\textbf{.3}}&-.3&\textcolor{blue}{0}&.9&\textcolor{blue}{.9}&.1&\textcolor{blue}{1}\\
& \textcolor{blue}{1}&&\textcolor{blue}{.7}&&\textcolor{blue}{.4}&&\textcolor{blue}{.9}&&\textcolor{blue}{1}&\\
\textcolor{blue}{0}&0&\textcolor{blue}{0}&.3&\textcolor{blue}{.3}&.6&\textcolor{blue}{.9}&.1&\textcolor{blue}{1}&0&\textcolor{blue}{1}\\
& \textcolor{blue}{1}&&\textcolor{blue}{1}&&\textcolor{blue}{1}&&\textcolor{blue}{1}&&\textcolor{blue}{1}&
\end{array}
\right)
\end{array}$
\caption[A circuit in a matrix in $P(n)$]{A matrix in $P(n)$ along with the matrix rewritten with the partial sums between the entries and a circuit of inner partial sums shown in boldface red}
\label{fig:asmcircuit}
\end{figure}

Label the corner matrix entries in the circuit alternately ($+$) and ($-$). Define 
\[k'=\min(r_{ij},1-r_{i'j'},c_{i''j''},1-c_{i''',j'''})\]
where $r_{ij}$, $r_{i'j'}$, $c_{i''j''}$, and $c_{i''',j'''}$ are taken over respectively the row partial sums to the right of a ($-$) corner along the circuit, the row partial sums to the right of a ($+$) corner, the column partial sums below a ($-$) corner, and the column partial sums below a ($+$) corner. Subtract $k'$ from the entries labeled ($-$) and add $k'$ to
the entries labeled ($+$). Subtracting and adding $k'$ in this way
preserves the row and column sums and keeps all the partial sums weakly between 0 and 1 (satisfying (\ref{eq:partialcolumnsum})--(\ref{eq:sumagain2})), so the result is another matrix $X'$ in $P(n)$ with at 
least one more non--inner partial sum than $X$. 

Now give opposite labels to the corners in the circuit in $X$ and subtract and add another constant $k''$ in a similar way to obtain another matrix $X''$ in $P(n)$ with at least one more non--inner partial sum than $X$. Then $X$ is a convex combination of $X'$ and $X''$, namely
$X=\frac{k''}{k'+k''} X' + \frac{k'}{k'+k''}X''$.
Therefore, by repeatedly applying this procedure, $X$ can be written as a convex combination of alternating sign matrices (i.e. matrices of $P(n)$ with no inner partial sums). 
\end{proof}

\section{Properties of the ASM polytope}
\label{sec:properties}

Now that we can describe the alternating sign matrix polytope in terms of inequalities, let us use this inequality description to examine some of the properties of $ASM_n$, namely, its facets, its vertices, and its projection to the permutohedron.

To make the proofs of the next two theorems more transparent, we introduce modified square ice configurations called \emph{simple flow grids} which will be used more extensively in Section~\ref{sec:facelattice}. 
Consider a directed graph with $n^2+4n$ vertices: $n^2$ `internal' vertices $(i,j)$ and $4n$ `boundary' vertices $(i,0)$, $(0,j)$, $(i,n+1)$, and $(n+1,j)$ where $i,j=1,\ldots,n$. These vertices are naturally depicted in a grid in which vertex $(i,j)$ appears in row $i$ and column $j$. Define the \emph{complete flow grid} $C_n$ to be the directed graph on these vertices with edge set $\{((i,j),(i,j\pm1)),((i,j),(i\pm1,j))\}$ for $i,j=1,\ldots,n$. So $C_n$ has directed edges pointing in both direction between neighboring internal vertices in the grid, and also directed edges from internal vertices to neighboring border vertices.

\begin{definition}
\label{definition:simpleflowgrid}
A \emph{simple flow grid} of order $n$ is a subgraph of $C_n$ consisting of all the vertices of $C_n$ for which four edges are incident to each internal vertex: either four edges
directed inward, four edges directed outward, or two horizontal edges pointing in the
same direction and two vertical edges pointing in the same direction.
\end{definition}

\begin{figure}[htbp]
\centering
\[
\includegraphics[scale=0.2]{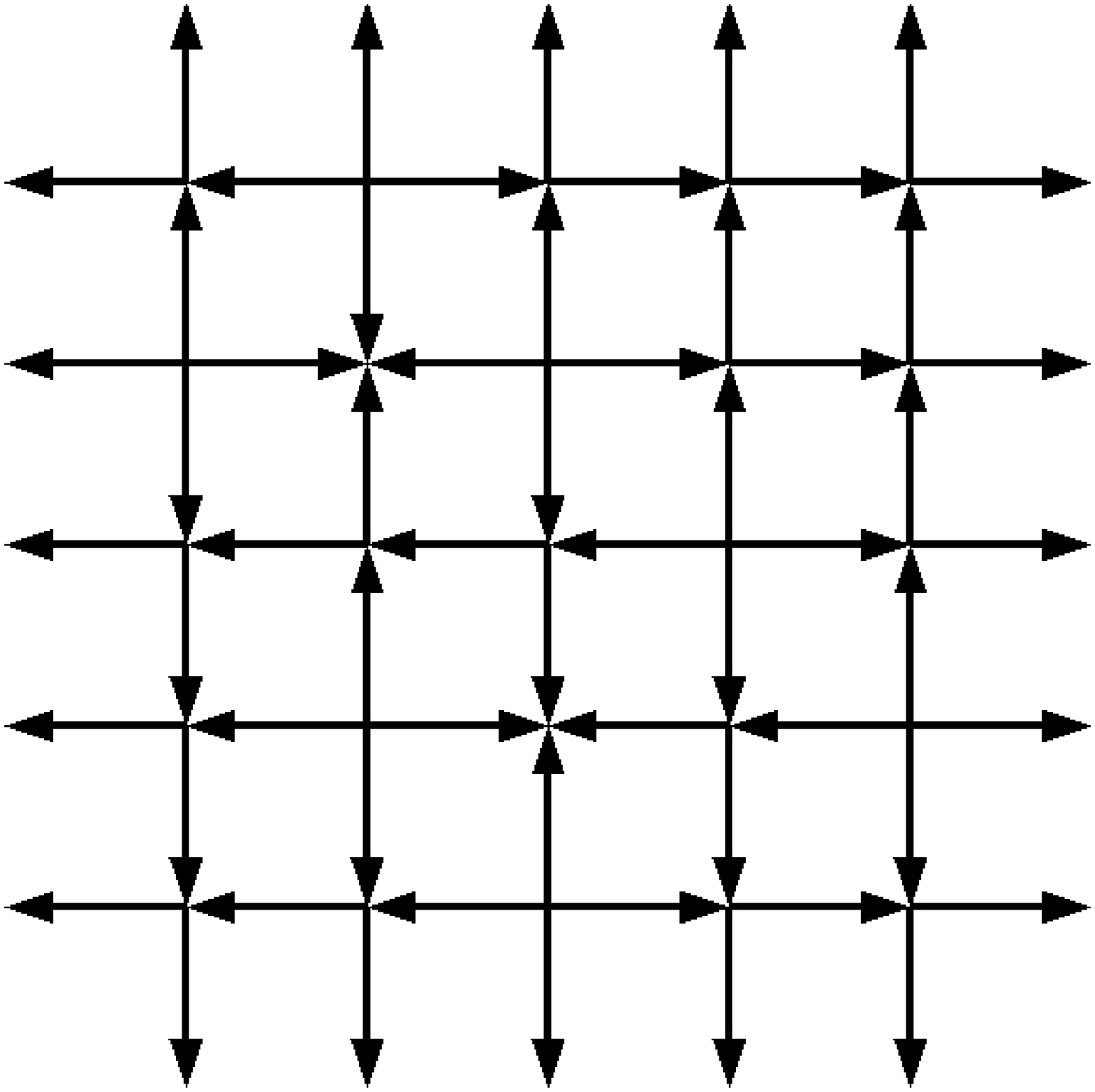}
\hspace{.7in}
\includegraphics[scale=0.25]{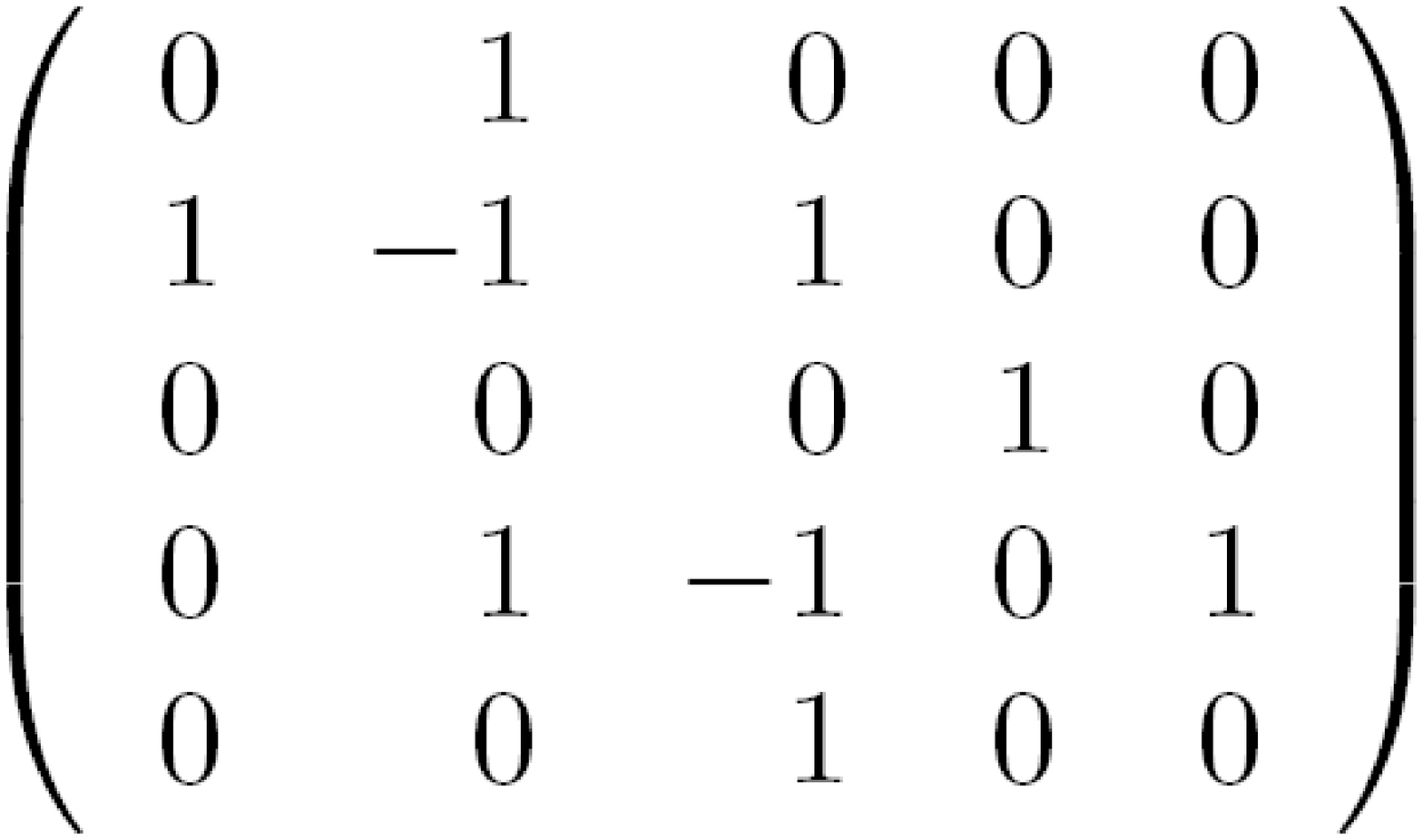}
\]
\caption{The simple flow grid---ASM correspondence}
\label{fig:5by5grid}
\end{figure}

\begin{proposition}
\label{prop:sfgasmbij}
There exists an explicit bijection between simple flow grids of order $n$ and $n\times n$ alternating sign matrices.
\end{proposition}
\begin{proof}
Given an ASM $A$, we will define a corresponding directed graph $g(A)$ on the $n^2$ internal vertices and $4n$ boundary vertices arranged on a grid as described above. Let each entry $a_{ij}$ of $A$ correspond to the internal vertex $(i,j)$ of $g(A)$.
For neighboring vertices $v$ and $w$ in $g(A)$ let there be a directed edge from $v$ to $w$ if the partial sum from the border of the matrix to the entry corresponding to $v$ in the direction pointing toward $w$ equals 1. By the definition of alternating sign matrices, there will be exactly one directed edge between each pair of neighboring internal vertices and also a directed edge from an internal vertex to each neighboring border vertex. Vertices of $g(A)$ corresponding to 1's are sources and vertices corresponding to $-1$'s are sinks. The directions of the rest of the edges in $g(A)$ are determined by the placement of the 1's and $-1$'s, in that there is a series of directed edges emanating from the 1's and continuing until they reach a sink or a border vertex. Thus $g(A)$ is a simple flow grid. Also, given a simple flow grid we can easily find the corresponding ASM by replacing all the sources with 1's and all the sinks with $-1$'s. Thus simple flow grids are in one-to-one correspondence with ASMs (see Figure~\ref{fig:5by5grid}).
\end{proof}

Simple flow grids are, in fact, almost the same as configurations of the six-vertex model of square ice with domain wall boundary conditions (see the discussion in~\cite{BRESSOUDBOOK}), the only difference being that the horizontal arrows point in the opposite direction. 

Recall that for $n\ge 3$ the Birkhoff polytope has $n^2$ facets (faces of dimension one less than the polytope itself). ($B_2=ASM_2$ is simply a line segment, so the number of facets equals the number of vertices which is 2.) Each facet of the Birkhoff polytope consists of all nonnegative doubly stochastic matrices with a zero in a fixed entry, that is, where one of the defining inequalities is made into an equality. The analogous theorem for $ASM_n$ is the following.

\begin{theorem}
\label{theorem:facets}
$ASM_n$ has $4[(n-2)^2+1]$ facets, for $n\ge 3$.
\end{theorem}

\begin{proof}
Note that the $4 n^2$ defining inequalities for $X\in ASM_n$ given in (\ref{eq:partialcolumnsum}) and (\ref{eq:partialrowsum}) can be restated as
\[\sum_{i'=1}^{i} x_{i'j} \ge 0 \hspace{.5 in}\sum_{j'=1}^{j} x_{ij'} \ge 0\]
\[\sum_{i'=i}^n x_{i'j} \ge 0 \hspace{.5 in}\sum_{j'=j}^n x_{ij'} \ge 0\]
for $i,j=1,\ldots,n$.
We have rewritten the statement that the row and column partial sums from the left or top must be less than or equal to $1$ as the row and column partial sums from the right and bottom must be greater than or equal to $0$. By counting these defining inequalities, one sees that there could be at 
most $4n^2$ facets, each determined by making one of the above inequalities an equality. It is left to determine how many of these equalities 
determine a face of dimension less than $(n-1)^2-1$. 

By symmetry we can determine the number of facets coming from the inequalities $\sum_{i'=1}^{i} x_{i'j} \ge 0$ for $i,j=1,\ldots n$ and then multiply by 4. Since the full row and column sums always equal $1$, the equalities such as $\sum_{i'=1}^n x_{i'j} = 0$ yield the empty face ($i=n$).
Also, $\sum_{i'=1}^{n-1} x_{i'j}\ge 0$ is implied from the fact that $x_{nj}\ge0$ ($i=n-1$). The inequalities $x_{i'1}\ge 0$ for all $i'$, i.e. the entries in the first column are nonnegative, imply that $\sum_{i'=1}^i x_{i'1}\ge 0$ for $2\le i\le n-1$ ($j=1$), thus each of these sets is a face of dimension less than $(n-1)^2-1$, and similarly for $\sum_{i'=1}^i x_{i'n}\ge 0$ for $2\le i\le n-1$ ($j=n$) the partial sums of the last column.

So we are left with the $(n-2)^2$ inequalities $\sum_{i'=1}^{i} x_{i'j} = 0$ for $i=1,\ldots n-2$ and $j=2,\ldots,n$ along with the inequality $x_{11}\ge 0$. For our symmetry argument to work, we do not include $x_{n1}\ge 0$ in our count since $x_{n1}\ge 0$ is also an inequality of the form $\sum_{j'=1}^{j} x_{ij'} \ge 0$.

Thus $ASM_n$ has at most $4[(n-2)^2+1]$ facets, given explicitly by the $4(n-2)^2+4$ sets of all $X\in ASM_n$ which satisfy one of the following:
\begin{equation}
\label{eq:faceteq}
%\sum_{i'=1}^{i} x_{i'j} &=& 0, \hspace{.2 in}\sum_{j'=1}^{j} x_{ij'} = 0, \hspace{.2 in}i\in\{1,\ldots,n-2\},j\in\{2,\ldots,n-1}\\
%\sum_{i'=i}^n x_{i'j} &=& 0, \hspace{.2 in}\sum_{j'=j}^n x_{ij'} = 0,\hspace{.2 in}  i,j\in\{3,\ldots,n\},\\
\sum_{i'=1}^{i-1} x_{i'j}=0, \sum_{j'=1}^{j-1} x_{ij'}=0,
\sum_{i'=i+1}^n x_{i'j}=0, \sum_{j'=j+1}^n x_{ij'}=0, i,j\in\{2,\ldots,n-1\},
\end{equation}
\begin{equation}
\label{eq:corners}
x_{11}=0,\mbox{ } x_{1n}=0,\mbox{ } x_{n1}=0,\mbox{ or } x_{nn}=0.
\end{equation}

They are facets (not just faces) since each equality determines exactly one more entry of the matrix, decreasing the dimension by one.

Recall that a directed edge in a simple flow grid $g(A)$ represents a location in the corresponding ASM $A$ where the partial sum equals 1, thus a directed edge missing from $g(A)$ represents a location in $A$ where the partial sum equals 0. Thus we can represent each of the $4(n-2)^2$ facets of (\ref{eq:faceteq}) as subgraphs of the complete flow grid $C_n$ from which a single directed edge has been removed: $((i\pm 1,j),(i,j))$ or $((i,j\pm 1),(i,j))$ with $i,j\in\{2,\ldots,n-1\}$. We can represent the facets of (\ref{eq:corners}) as subgraphs of $C_n$ from which two directed edges have been removed: $((1,1),(1,2))$ and $((1,1),(2,1))$, $((1,n),(1,n-1))$ and $((1,n),(2,n))$, $((n,1),(n-1,1))$ and $((n,1),(n,2))$, or $((n,n),(n,n-1))$ and $((n,n),(n-1,n))$.

Now given any two facets $F_1$ and $F_2$, it is easy to exhibit a pair of ASMs $\{X_1,X_2\}$ such 
that $X_1$ lies on $F_1$ and not on $F_2$. 
Include the directed edge(s) corresponding to $F_2$ but not the directed edge(s) corresponding to $F_1$ in $g(X_1)$, then do the opposite for $X_2$. Thus each of the 
$4[(n-2)^2+1]$ equalities gives rise to a unique facet.
\end{proof}

\begin{corollary}
For $n\ge 3$, the number of facets of $ASM_n$ on which an ASM $A$ lies is given by
$2(n-1)(n-2)+(\mbox{number of corner 1's in $A$}).$
\end{corollary}
\begin{proof}
Each $0$ around the border of $A$ represents one facet. Thus the number of facets corresponding to border zeros of $A$ equals $4(n-1) - (\mbox{\# 1's around the border of $A$})$. Then there are $2(n-2)(n-3)$ facets represented by directed edges pointing in the opposite directions to the directed edges in the $(n-2)\times(n-2)$ interior array of $g(A)$. The sum of these numbers gives the above count.
\end{proof}

Even though $ASM_n$ is defined as the convex hull of the ASMs, it requires some proof that each ASM is actually an extreme point of $ASM_n$. 
\begin{theorem}
The vertices of $ASM_n$ are the $n\times n$ alternating sign matrices.
\end{theorem}

\begin{proof} 
Fix an $n\times n$ ASM $A$. In order to show 
that $A$ is a vertex of $ASM_n$, we need to find a 
hyperplane with $A$ on one side and all the other ASMs on the other side. Then since $ASM_n$ is the convex hull of $n\times n$ ASMs, $A$ would necessarily be a vertex. 

Consider the simple flow grid corresponding to $A$. 
 In any simple flow grid there are, by definition, $2n(n+1)$ directed edges, where for each entry of the corresponding ASM there is a directed edge whenever the partial sum in that direction up to that point equals 1. 
Since the total number of directed edges in a simple flow grid is fixed, $A$ is the only ASM with all of those partial sums equaling 1. Thus the hyperplane where the sum of those partial sums 
equals $2n(n+1)-\frac{1}{2}$ will 
have $A$ on one side and all the other ASMs on the other. Thus the $n\times n$ ASMs 
are the vertices of $ASM_n$.
\end{proof}

Another interesting property of the ASM polytope is its relationship to the 
permutohedron. For a vector $z = (z_1, z_2, \ldots, 
z_n)\in\mathbb{R}^n$ with distinct entries, define the permutohedron $P_z$ as the convex hull of all vectors obtained by permuting the entries of $z$. That is,
\begin{equation}
\label{eq:permut}
P_z = \mbox{conv}\{(z_{\omega (1)}, z_{\omega (2)}, \ldots, z_{\omega (n)})\mbox{ }|\mbox{ } \omega \in S_n\}.
\end{equation}
Also, for such a vector $z$, let $\phi_z$ be the mapping from the set of $n\times n$ real matrices to $\mathbb{R}^n$ defined by
\[\phi_z (X)= z X, \mbox{for any $n\times n$ real matrix $X$.}\]

It is well known, and follows immediately from the definitions, that $P_z$ is the image of the Birkhoff polytope under the projection $\phi_z$. 
\begin{proposition}
\label{prop:birkpermut}
Let $B_n$ be the Birkhoff polytope and $z$ be a vector in $\mathbb{R}^n$ with distinct entries. Then 
\begin{equation}
\phi_z(B_n)=P_z.
\end{equation}
\end{proposition}
This result is one of many classical results about the Birkhoff polytope dating back to Hardy, Littlewood, and P\'olya~\cite{HARDYLP1}~\cite{HARDYLP2}.  See~\cite{MIRSKY} for a nice summary of relevant results.

The next theorem states that when the same projection map is applied to 
$ASM_n$, the image is the same permutohedron whenever $z$ is a decreasing vector.
For the proof of this theorem we will need the concept of majorization~\cite{MAJORIZATION}. 
\begin{definition}
\label{def:maj}
Let $u$ and $v$ be vectors of length $n$. 
Then $u\preceq v$ (that is $u$ is \emph{majorized} by $v$) if
\begin{equation}
\begin{cases}
\sum_{i=1}^k u_{[i]} \le \sum_{i=1}^k v_{[i]},&\mbox{for } 1\le k\le n-1\\
\sum_{i=1}^n u_i = \sum_{i=1}^n v_i&\end{cases}
\end{equation}
where the vector $(u_{[1]},u_{[2]},\ldots,u_{[n]})$ is obtained from $u$ by rearranging its components so that they are in decreasing order, and similarly for $v$.
\end{definition}

\begin{theorem} 
Let $z$ be a decreasing vector in $\mathbb{R}^n$ with distinct entries. 
Then 
\begin{equation}
\phi_z(ASM_n)=P_z.
\end{equation}
\end{theorem}
\begin{proof}
It follows from Proposition~\ref{prop:birkpermut} and $B_n\subseteq ASM_n$ that $P_z\subseteq \phi_z(ASM_n)$. Thus it only remains to be shown that $\phi_z(ASM_n)\subseteq P_z$.

Let $z$ be a  decreasing $n$--vector (so that $z_i = z_{[i]}$) and $X = \{x_{ij}\}$ an $n\times n$ ASM. Then there is a proposition of Rado which states that for vectors $u$ and $v$ of length $n$, $u\preceq v$ if and only if $u$ lies in the convex hull of the $n!$ permutations of the entries of 
$v$~\cite{RADO}.
Therefore the proof will be completed by showing $z X\preceq z$. By Definition~\ref{def:maj} we need to show
\begin{align}
\label{eq:proj1}
\sum_{j=1}^k (z X)_{[j]}&\le \sum_{j=1}^k z_{j}, \hspace{.4in} 1\le k\le n-1\\
\label{eq:proj2}
\sum_{j=1}^n (z X)_j&= \sum_{j=1}^n z_j
\end{align}
where the $j$th component $(z X)_j$ of $z X$ is given by $\sum_{i=1}^n z_i x_{ij}$. 

To verify (\ref{eq:proj2}) note that since $\sum_{j=1}^n x_{ij}=1$,
\[\sum_{j=1}^n (z X)_j = \sum_{j=1}^n \sum_{i=1}^n z_i x_{ij} = \sum_{i=1}^n z_i \sum_{j=1}^n x_{ij} = \sum_{i=1}^n z_i.\]

To prove (\ref{eq:proj1}) we will show that $\sum_{j\in J}  (z X)_{j} \le \sum_{j=1}^{|J|} z_j$ given any $J\subseteq \{1,\ldots, n\}$, so that in particular $\sum_{j=1}^{|J|} (z X)_{[j]}\le \sum_{j=1}^{|J|} z_j$. 

We will need to verify the following:
\begin{align}
\label{eq:vim}
\sum_{i=1}^m \sum_{j\in J} x_{ij}&\le \min(m,|J|) \hspace{.4in} \forall \mbox{ }m\in \{1,\ldots, n\}.\\
\label{eq:vik}
\sum_{i=1}^n \sum_{j\in J} x_{ij}&= |J|.
\end{align}

To prove (\ref{eq:vim}) note that 
\[
\sum_{i=1}^m \sum_{j\in J} x_{ij} =  \sum_{j\in J} \sum_{i=1}^m x_{ij} \le |J|
\]
since $\sum_{i=1}^m x_{ij}\le 1$. But also, since $\sum_{i=1}^m x_{ij}\ge 0$ and $\sum_{j=1}^n  x_{ij} = 1$ we have that
\[
\sum_{j\in J} \sum_{i=1}^m x_{ij} \le \sum_{j=1}^n \sum_{i=1}^m x_{ij} = \sum_{i=1}^m \sum_{j=1}^n  x_{ij} = m.
\]
To prove (\ref{eq:vik}) observe, 
\[
\sum_{i=1}^n \sum_{j\in J} x_{ij} = \sum_{j\in J} \sum_{i=1}^n x_{ij} = \sum_{j\in J} 1 = |J|
\] 
since the columns of $X$ sum to 1. 
Therefore using (\ref{eq:vim}) and~(\ref{eq:vik}) we see that

\begin{align*}
\sum_{i=1}^n z_i x_{ij}&= \sum_{j\in J} \sum_{i=1}^n z_i x_{ij} =  \sum_{i=1}^n z_i \sum_{j\in J} x_{ij} = \sum_{k=1}^{n-1} (z_k-z_{k+1}) \sum_{i=1}^k \sum_{j\in J} x_{ij} + z_n \sum_{i=1}^n \sum_{j\in J} x_{ij}\\
&= \sum_{k=1}^{n-1} (z_k-z_{k+1}) \sum_{i=1}^k \sum_{j\in J} x_{ij} + z_n |J| && 
\text{by~(\ref{eq:vik})}\\
&= \sum_{k=1}^{|J|-1} (z_k-z_{k+1}) \sum_{i=1}^k \sum_{j\in J} x_{ij} + \sum_{k=|J|}^{n-1} (z_k-z_{k+1}) \sum_{i=1}^k \sum_{j\in J} x_{ij} + z_n |J|\\
&\le \sum_{k=1}^{|J|-1} (z_k-z_{k+1}) k + \sum_{k=|J|}^{n-1} (z_k-z_{k+1}) |J| + z_n |J| &&\text{by~(\ref{eq:vim})}\\
&\le \sum_{k=1}^{|J|-1} z_k - z_{|J|} (|J|-1) + (z_{|J|}-z_{n})|J| + z_n |J|\\
&\le \sum_{k=1}^{|J|} z_k.
\end{align*}

Thus $z X \preceq z$ and so $z X$ is contained in the convex 
hull of the permutations of $z$. Therefore $\phi_z(ASM_n)=P_z$.
\end{proof}

\section{The face lattice of the ASM polytope}
\label{sec:facelattice}

Another nice result about the Birkhoff polytope is the structure of its face lattice~\cite{BILLERA}. 
Associate to each permutation matrix $X$ a bipartite graph with vertices $u_1,u_2,\ldots,u_n$ and $v_1,v_2,\ldots,v_n$ where there is an edge connecting $u_i$ and $v_j$ if and only if there is a 1 in the $(i,j)$ position of $X$. Such a graph will be a perfect matching on the complete bipartite graph $K_{n,n}$. A graph $G$ is called elementary if every edge is a member of some perfect matching of~$G$.
\begin{theorem}[Billera--Sarangarajan]
\label{prop:elementary}
The face lattice of the Birkhoff polytope is isomorphic to the lattice of 
elementary subgraphs of $K_{n,n}$ ordered by inclusion.
\end{theorem}
This lattice structure was first identified by Billera and Sarangarajan in~\cite{BILLERA} and~\cite{BILLERAFPSAC}, but the set of faces itself
was first characterized and studied extensively 
by Brualdi and Gibson in~\cite{BRUALDIGIB1} and~\cite{BRUALDIGIB2}
using certain 0-1 matrices which correspond trivially to elementary subgraphs of $K_{n,n}$.
Other relevant results were also obtained by Balinski and Russakoff in~\cite{BALRUS1} and~\cite{BALRUS2}.

A similar statement can be made about the face lattice of the ASM polytope using simple flow grids (see Definition~\ref{definition:simpleflowgrid}) in place of perfect matchings, the complete flow grid $C_n$ instead of the complete bipartite graph $K_{n,n}$, and \emph{elementary flow grids} in place of elementary graphs.

\begin{definition}
An \emph{elementary flow grid} $G$ is a subgraph of the complete flow grid $C_n$ such that the edge set of $G$ is the union of the edge sets of simple flow grids.
\end{definition}

Now for any face $F$ of $ASM_n$ define the grid corresponding to the 
face, $g(F)$, to be the union over all the vertices of $F$ of the simple flow grids corresponding to the vertices. That is,
\[g(F) = \bigcup_{vertices\mbox{ }A\in F} g(A).
\]
Thus $g(F)$ is an elementary flow grid since its edge set is the union of the edge sets of simple flow grids.

Now we wish to define the converse, that is, given an elementary flow grid $G$ 
we would like to know the corresponding face $f(G)$ of $ASM_n$. Define $f(G)$ to be the convex hull of the vertices of $ASM_n$ whose corresponding simple flow grids are contained in the elementary flow grid $G$. So let
\[f(G) = \mbox{conv}\{\mbox{vertices }A\in ASM_n \mbox{ }|\mbox{ } g(A)\subseteq G\}.\]

Recall that we can represent each of the facets of $ASM_n$ either as subgraphs of the complete flow grid $C_n$ from which one of the directed edges in the set $S=\{((i\pm 1,j),(i,j))$, $((i,j\pm 1),(i,j))$ $|$ $i,j\in\{2,\ldots,n-1\}\}$ has been removed or from which one of the pairs of directed edges in the set $T=\{\{((1,1)(1,2)),((1,1),(2,1))\}$, $\{((1,n),(1,n-1)),((1,n),(2,n))\}$, $\{((n,1),(n-1,1)),((n,1),(n,2))\}$, $\{((n,n),(n,n-1)),((n,n),(n-1,n))\}\}$ has been removed.

Thus each of the directed edges in $S$ and the first of each pair of directed edges in $T$ that are not in $G$ represent facets that contain $f(G)$. Let the collection of these directed edges be called $\{e_1, e_2, \ldots, e_k\}$ and their corresponding facets $\{F_1, F_2, \ldots, F_k\}$. Let $I=\bigcap_{j=1}^k F_j$ be the intersection of these facets. Thus $I$ is a face of $ASM_n$ and $f(G)\subseteq I$. 

We wish to show that $f(G)$ equals $I$. So suppose $f(G)\subsetneq I$. Then since $I$ is a face of $ASM_n$ and $f(G)$ is defined as the convex hull of vertices of $ASM_n$ there exists an additional vertex $B\in I$ of $ASM_n$ such that $B\notin f(G)$. But $g(B)$ must be missing the directed edges $e_1,e_2,\ldots,e_k$ since $B\in I$, thus all the directed edges of $g(B)$ must be in $G$. Therefore $g(B)\subseteq G$ so that $B\in f(G)$ which is a contradiction. So $f(G)=I$.
Thus $f(G)$ is a face of $ASM_n$ since it is the intersection of faces of $ASM_n$. 

It can easily be seen that $f(g(F))=F$ and $g(f(G))=G$. Also if $F_1$ and $F_2$ are faces of $ASM_n$ then $F_1\subseteq F_2$ if and only if $g(F_1)\subseteq g(F_2)$.

Thus elementary flow grids are in bijection with the faces of $ASM_n$ (if we also regard the empty grid as an elementary flow grid). Elementary flow grids can be made into a lattice by 
inclusion, where the join is the union of the edge sets and the meet is 
the largest elementary flow grid made up of the directed edges from the 
intersection of the edges sets.

This discussion yields the following theorem:
\begin{theorem}
\label{thm:facelattice}
The face lattice of $ASM_n$ is isomorphic to the 
lattice of all $n\times n$ elementary flow grids (or equivalently all $n\times n$ square ice configurations with domain wall boundary conditions) ordered by inclusion.
\end{theorem}

The dimension of any face of $ASM_n$ can be determined by looking at $g(F)$ as in the following theorem. 
The characterization of edges of $ASM_n$ is analogous to the result for the Birkhoff polytope which states that the graphs representing edges of $B_n$ are the elementary subgraphs of $K_{n,n}$ which have exactly one cycle~\cite{BALRUS2}~\cite{BILLERA}~\cite{BRUALDIGIB2}.

\begin{figure}[htp]
\centering
\includegraphics[scale=0.25]{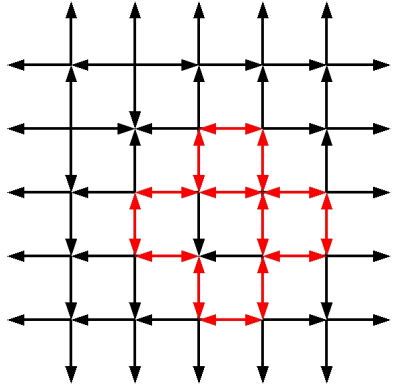}
\caption[An elementary flow grid containing 3 doubly directed regions]{An elementary flow grid containing 3 doubly directed regions which corresponds by Theorem~\ref{thm:dim} to a 3-dimensional face of $ASM_5$}
\label{fig:redbluecircuits}
\end{figure}

Given an elementary flow grid $G$, define a \emph{doubly directed region} as a collection of cells in $G$ 
completely bounded by double directed edges but containing no double directed edges in the interior (see Figure~\ref{fig:redbluecircuits}). Let $\alpha(G)$ denote the number of doubly directed regions in $G$.

\begin{theorem}
\label{thm:dim}
The dimension of a face $F$ of $ASM_n$ is the number of doubly directed regions in the corresponding elementary flow grid $g(F)$.
In particular, the edges of $ASM_n$ are represented by elementary flow grids containing exactly one cycle of double directed edges. 
\end{theorem}

\begin{proof}
We proceed by induction on the dimension of the face of $ASM_n$. The simple flow grid corresponding to any ASM $A$ has no double directed edges, thus $\alpha(g(A))=0$. Now
suppose for every $m$--dimensional face of $ASM_n$, the number of doubly directed regions of the elementary flow grid corresponding to the face equals $m$. 
Let $F$ be an $(m+1)$--dimensional face of $ASM_n$ and $F'$ an $m$--dimensional subface of $F$. 
We assume $\alpha(g(F')) = m$ and wish to show that $\alpha(g(F)) = m+1$. 

Now $g(F)$ is the elementary flow grid whose edge set is the union of the edge sets of $g(F')$ and $g(A)$ over all ASMs $A$ in $F-F'$. 
Every vertex in a simple flow grid must have even indegree and even outdegree. Therefore, if we wish to obtain $g(A)$ from $g(A')$, where $A'$ is an ASM in $F'$, by reversing some directed edges, 
the number of directed edges reversed at each vertex must be even. Thus taking the union of the directed edges of $g(A')$ with the directed edges of 
$g(A)$ forms one or more circuits of double directed edges, where at least one of the double directed edges is not in $g(F')$.
Therefore $g(F)$ has at least one more doubly directed region than $g(F')$,
so $\alpha(g(F))\ge m+1$.
Then since $g(ASM_n)$ equals the complete flow grid $C_n$, we have that 
$\alpha(g(ASM_n)) =\alpha(C_n)= (n-1)^2 = \mbox{dim}(ASM_n)$.
Therefore moving up the face lattice one rank increases the number of doubly directed regions by exactly one, so $\alpha(g(F)) = m+1$.
\end{proof}

\begin{figure}[htp]
\centering
\includegraphics[scale=0.25]{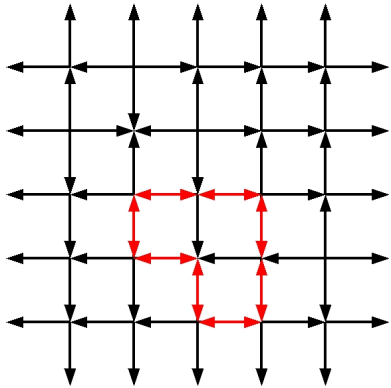}
\caption{The elementary flow grid representing an edge in $ASM_5$}
\label{fig:redcircuit}
\end{figure}

See Figure~\ref{fig:redcircuit} for the elementary flow grid representing 
the edge in $ASM_5$ between
\[
\left(
\begin{array}{rrrrr}
0 & 1 & 0 & 0 & 0\\
1 & -1 & 1 & 0 & 0\\
0&0&0&1&0\\
0&1&-1&0&1\\
0&0&1&0&0
\end{array}\right)
\mbox{ and }
\left(
\begin{array}{rrrrr}
0 & 1 & 0 & 0 & 0\\
1 & -1 & 1 & 0 & 0\\
0&1&0&0&0\\
0&0&0&0&1\\
0&0&0&1&0
\end{array}
\right).
\]

\section{Acknowledgments}
This work is based on research which is a part of the author's doctoral thesis at the University of Minnesota under the direction of Dennis Stanton. The author would like to thank Professor Stanton for the many helpful discussions and encouragement. 
%The author is also very grateful to the referee for the detailed report which helped to improve the clarity of the manuscript.

%-------------------------------------------BIBLIOGRAPHY------------------------------------------------------------%

\bibliographystyle{amsalpha}

\end{document}